\documentclass[12pt]{amsart}

\usepackage{amsthm,amssymb,amsmath,amsfonts}
\usepackage{enumitem}
\usepackage{color}
\usepackage{graphicx}
\usepackage{caption}
\usepackage{subcaption}

\newtheorem{theorem}{Theorem}
\newtheorem{proposition}{Proposition}
\newtheorem{corollary}{Corollary}

\title{The Hamiltonian Dynamics of Planar Magnetic Confinement}
\author{Gabriel Martins}

\begin{document}
\maketitle


\begin{abstract}
Inspired by a question of Colin de Verdi\`{e}re and Truc we study the dynamics of a classical charged particle moving in a bounded planar domain $\Omega$ under the influence of a magnetic field $\mathbf{B}$ which blows up at the boundary of the domain. We prove that under appropriate blow-up conditions the particle will never reach the boundary. As a corollary we obtain completeness of the magnetic flow. Our blow-up condition is that $\mathbf{B}$ should not be integrable along normal rays to the boundary, while its tangential derivative should be integrable along those same rays.
\end{abstract}

\tableofcontents


\section{Introduction}
\label{sec:intro}
We study the global asymptotics of a charged  particle in a planar region under the influence of a magnetic field. Much work has been done on this and related subjects, see \cite{Arnold,Arnold2,Braun,Castilho,deVerdiere,Montgomery,Truc}, but our main motivation comes from the work \cite{deVerdiere}, where Colin de Verdi\`{e}re and Truc establish growth conditions for a magnetic field $\mathbf{B}$ defined on a Riemannian domain $\Omega$ with boundary satisfying certain regularity assumptions, so that the associated magnetic Schr\"{o}dinger operator is essentially self-adjoint on $C^\infty_0(\Omega)$, which means any quantum particle supported in $\Omega$ is confined to $\Omega$ for all time (for more details see section \ref{sec:quantum}). In that same work they ask if the analogous classical statement is true. We give here a positive answer to their question in dimension 2.

The results obtained here can be viewed as a classical 2 dimensional toy model for the dynamics inside of a Tokamak, a magnetic confinement device used in fusion power generators. Our approximation excludes collision of the interior particles with the Tokamak (the boundary of our region). We also have ongoing work on higher dimensional analogues of our results.

A magnetic field on a domain $\Omega \subset \mathbb{R}^n$ is given by a 2-form $\mathbf{B}$ defined over $\Omega$. The 2-form is equivalent to an antisymmetric endomorphism $Y:T\mathbb{R}^n\to T\mathbb{R}^n$ via the correspondence $\mathbf{B}_q= \langle\cdot,Y_q(\cdot)\rangle$. The corresponding equations of motion for a particle of charge $e$ and mass $m$ moving in $\Omega$ under the influence of $\mathbf{B}$ are:
  \begin{equation}\label{eq:genmag}
  m\ddot{q} = eY_q(\dot{q})
  \end{equation}

If $\vec{B}$ is a vector field in $\mathbb{R}^3$ we may encode it into the 2-form $\mathbf{B}=\langle\cdot,\cdot\times\vec{B}\rangle$ then the endomorphism $Y$ above is simply $\cdot\times\vec{B}$ and equations (\ref{eq:genmag}) in this case take the familiar form:
  \begin{equation}
  m\ddot{q} = e\dot{q}\times\vec{B}
  \end{equation}
  
For the 2 dimensional picture one may consider a magnetic field of the form $\vec{B}=(0,0,B(x,y))$. This assumption forces particles in the $xy$-plane whose initial velocities are tangent to the plane to stay in this same plane for all time. The equations of a charged particle under the influence of this field are:
  \begin{equation}\label{eq:mag2d}
  m\ddot{q} = -eB(q)J\dot{q}
  \end{equation}
where $J = \left(\begin{smallmatrix}0&-1\\1&0\end{smallmatrix}\right)$ is the standard complex structure on $\mathbb{R}^2$. We will restrict our initial positions and our field to a bounded planar domain  $\Omega$  endowed with the Euclidian metric, and having smooth boundary $\partial\Omega$. If $\mathbf{B}(q) = B(q)dx \wedge dy$ is the magnetic field, then equations (\ref{eq:genmag}) of a charged particle moving in $\Omega$ under the influence of $\mathbf{B}$ are precisely equations (\ref{eq:mag2d}).


\section{Statement of main results}
\label{sec:state}
Choosing a potential 1-form $\mathbf{A}$, so that $\mathbf{B}=d\mathbf{A}$, we may describe the magnetic equations (\ref{eq:genmag}) as a Hamiltonian system on $T^*\Omega$ with its canonical symplectic structure and Hamiltonian:
  \begin{equation}\label{eq:Ham}
  H(q,p) = \frac{|p-e\mathbf{A}(q)|_*^2}{2m}
  \end{equation}
where $|\cdot|_*$ denotes the norm on $T^*\Omega$ induced by the Euclidian metric.

We will show that if the field $\mathbf{B}$ grows fast enough as it approaches the boundary, then a particle in $\Omega$ may never reach the boundary of the region in finite time. This implies in particular that the magnetic flow on $T^*\Omega$ is complete. In order to phrase our condition more precisely we introduce normal coordinates on a neighborhood of $\partial\Omega$.

Let $\mathcal{C}$ be a connected component of $\partial\Omega$, let $L$ be its length, let $\gamma:\mathbb{R}/L\mathbb{Z} \to \mathcal{C}$ be an arc length parametrization of this curve so that the normal vector $\nu(s)$ defined by $\ddot{\gamma}(s)=\kappa(s)\nu(s)$ points inward. We define normal coordinates $x:(0,N) \times \mathbb{R}/L\mathbb{Z} \to \Omega$ by
  \begin{equation}
  x(n,s) = \gamma(s) + n\nu(s)
  \end{equation}
    
For a small enough choice of $N>0$ the map $x$ is a diffeomorphism, we denote its image by $\Omega_\mathcal{C}(N) \subset \Omega$, which is a collar neighborhood for the boundary curve $\mathcal{C}$. We can now state our main theorem:


\begin{theorem}\label{main}
For every connected component $\mathcal{C}$ of $\partial\Omega$, write the magnetic field $\mathbf{B} = B(n,s) dn\wedge ds$ over $\Omega_\mathcal{C}(N)$ using normal coordinates. Suppose $\mathbf{B}$ is smooth and that for every $\mathcal{C}$ we have:
  \begin{equation} \label{eq:nintB}
  \lim_{n\to0} \left|\int_n^N B(m,s)dm\right| = \infty.
  \end{equation}
Furthermore suppose for all $\mathcal{C}$, there is a constant $D_{\mathcal{C}}$ such that:
  \begin{equation}\label{eq:spartial}
  \sup_s \int_0^N \left|\frac{\partial B}{\partial s}(m,s)\right| dm < D_{\mathcal{C}}
  \end{equation}
Then, no trajectory of the magnetic flow starting in $\Omega$ reaches the boundary $\partial\Omega$ in finite time.
\end{theorem}

Using this theorem we may obtain the simple corollary:

\begin{corollary}
The Hamiltonian flow defined by $H$ on $T^*\Omega$ is complete.
\end{corollary}

\textbf{Proof:} If a trajectory never approaches the boundary $\partial\Omega$, then it must be contained in a compact set, hence it must be defined for all time. Otherwise, it must approach the boundary, and by the previous theorem it must take infinite time to do so, in particular it must be defined for all time. $\square$\\

By restricting the form of the magnetic field we may give a more quantitative description of the boundary behavior of the charged particle. Denote the cotangent coordinates induced by the normal coordinates by $(n,s,p_n,p_s)$. In the next result we give an explicit lower bound for the distance a particle in $\Omega_\mathcal{C}(N)$ must keep from the boundary in any finite amount of time.


\begin{theorem}\label{quant}
Let $\mathcal{C}$ be a component of $\partial\Omega$, let $K=sup_s|\kappa(s)|$ be the maximum curvature of $\mathcal{C}$, let  $K'=sup_s|\kappa'(s)|$ and let $N>0$ be small so that we can define normal coordinates on $\Omega_\mathcal{C}(N)$ and such that $N<\varepsilon/K$, for some $0<\varepsilon<1$. Suppose the magnetic field has the form:
  \begin{equation}
  B(n,s) = \frac{M}{n^\alpha}+f(n,s), \quad \alpha\geq1
  \end{equation}
with $|f|\leq C_f$ a bounded smooth function with integrable $s$-partial
  \begin{equation}
  \sup_s\int_0^N\left| \frac{\partial f}{\partial s}(m,s) \right|dm \leq D_\mathcal{C}
  \end{equation}
  
Let $q(t)=(n(t),s(t))$ be a solution to the equations (\ref{eq:mag2d}) contained in $\Omega_\mathcal{C}(N)$ for $0\leq t\leq T$ with energy $H_0$. Then:
  \begin{equation}
  n(t) \geq  \begin{cases} (N^{-(\alpha-1)} + (\alpha-1)d(T))^{-\frac{1}{\alpha-1}}& \text{if } \alpha > 1\\
                                        Ne^{-d(T)}                         & \text{if }\alpha = 1 \end{cases}
  \end{equation}
where:
  \begin{equation}
  d(T) = D_0 + D_1T
  \end{equation}
with constants given by:
  \begin{equation}
  \begin{aligned}
  D_0 & = \frac{C_f}{|M|} + \frac{|p_s(0)|+\sqrt{2mH_0}(1+\varepsilon)}{e|M|}\\
  D_1 & = \sqrt{\frac{2H_0}{m}}\cdot\frac{D_\mathcal{C}}{|M|(1-\varepsilon)}+\frac{2H_0K'N}{e|M|(1-\varepsilon)}
  \end{aligned}
  \end{equation}
\end{theorem}


\section{Comparison to previous results}
\label{sec:comparison}
We now compare the conditions of theorem \ref{main} to the ones in \cite{deVerdiere}. To this purpose we will take $\Omega$ to be the unit disc for simplicity, we'll parametrize the boundary by $\gamma(s)=(\cos(s),\sin(s))$, let $r = \sqrt{x^2+y^2}$, $n = 1-r$ and set $e=m=1$. We use both cartesian and normal coordinates to express the magnetic field:
  \begin{equation}
  \mathbf{B} = B(x,y) dx\wedge dy = \tilde{B} dn\wedge ds
  \end{equation}
with this notation we have $\tilde{B}(n,s) = (n-1)B(x,y)$, since $dx\wedge dy = (n-1)dn\wedge ds$. The conditions we impose on the magnetic field in order to guarantee completeness are not stronger or weaker than the ones in \cite{deVerdiere}. One major difference is that we suppose the boundary of our domain to be smooth, while \cite{deVerdiere} works with more general domains. Aside from that, using the above notation, their condition amounts to:
  \begin{equation}
  |B(x,y)| \geq \frac{1}{n^2}
  \end{equation}

To begin the comparison with the hypothesis of theorem \ref{main}, we consider magnetic fields of the form:
  \begin{equation}
  \tilde{B}(n,s) = \frac{M}{n^\alpha}
  \end{equation}
in cartesian coordinates we have:
  \begin{equation}
  B(x,y) = \frac{M}{(n-1)n^\alpha} \sim \frac{M}{n^\alpha}
  \end{equation}
asymptotically as $n\to0$. For this type of magnetic field, while the conditions from \cite{deVerdiere} require $|M| \geq 1$ and $\alpha \geq 2$, our result is more flexible and only requires $M\neq 0$ and $\alpha \geq 1$.

On the other hand, the results from \cite{deVerdiere} don't require any control on the dependence of $\mathbf{B}$ on the variable $s$, while theorem \ref{main} assumes that the $s$-partial derivative of $\mathbf{B}$ is not too wild (i.e. condition (\ref{eq:spartial}) in theorem \ref{main}). For example the magnetic field
  \begin{equation}
  B(x,y) = \frac{2+\sin(s)}{n^2} = \frac{2n + y}{n^3}
  \end{equation}
satisfies their hypothesis but not ours, since $\partial B/\partial s$ is not integrable along normal rays. A very interesting problem is whether or not one might be able to remove the hypothesis (\ref{eq:spartial}) from theorem \ref{main}.


\section{Examples}\label{sec:ex}

Here are some examples of regions $\Omega$ and magnetic fields $B$ that satisfy our hypothesis. Again we denote $r=\sqrt{x^2+y^2}$.\\

\textbf{1.} We take $\Omega = \{q\in\mathbb{R}^2,|q|<R\}$, the disc of radius $R$ with a magnetic field
  \begin{equation}
  B(x,y) = \frac{M}{(R - r)^{\alpha}} + f(x,y)
  \end{equation}
with $f$ a smooth bounded function on the closed disc, $\alpha\geq1$ and $M\neq0$.\\

\textbf{2.} Take the domain to be an annulus $\Omega=\{q\in\mathbb{R}^2,R_1<|q|<R_2\}$ with magnetic field:
  \begin{equation}
  B(x,y) = \frac{M_2}{(R_2 - r)^{\alpha_2}} + \frac{M_1}{(r-R_1)^{\alpha_1}} + f(x,y)
  \end{equation}
with similar assumptions as before on $f$ and the constants $M_i$, $R_i$ and $\alpha_i$.

\begin{figure}
  \centering
  \resizebox{0.75\textwidth}{!}{
  \includegraphics{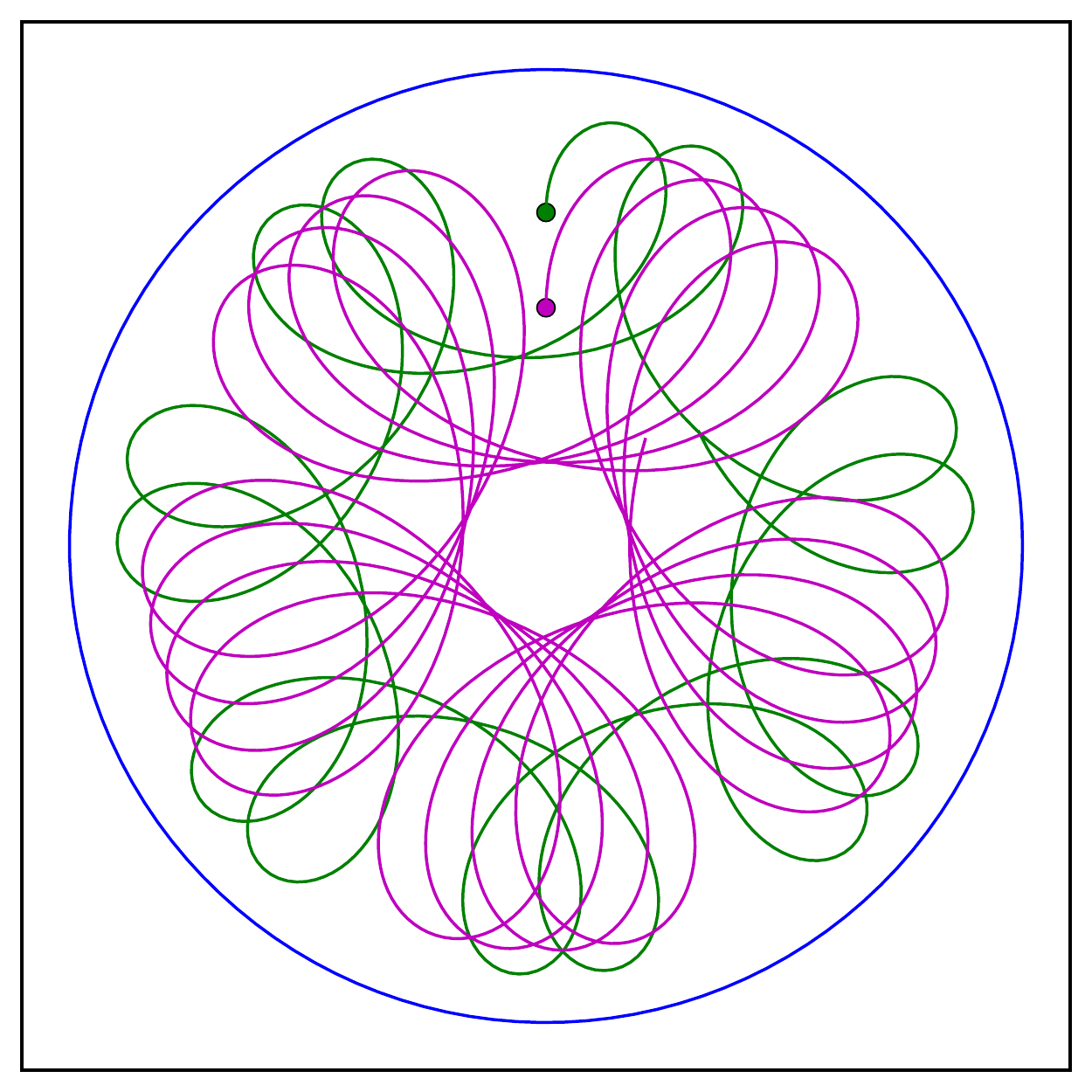}
  }
  \caption{$B(x,y)=\frac{1}{1-r}$}
  \label{fig:magsym}
\end{figure}

\begin{figure}
	\centering
	\begin{subfigure}[t]{0.45\textwidth}
		\centering
		\includegraphics[width=\textwidth]{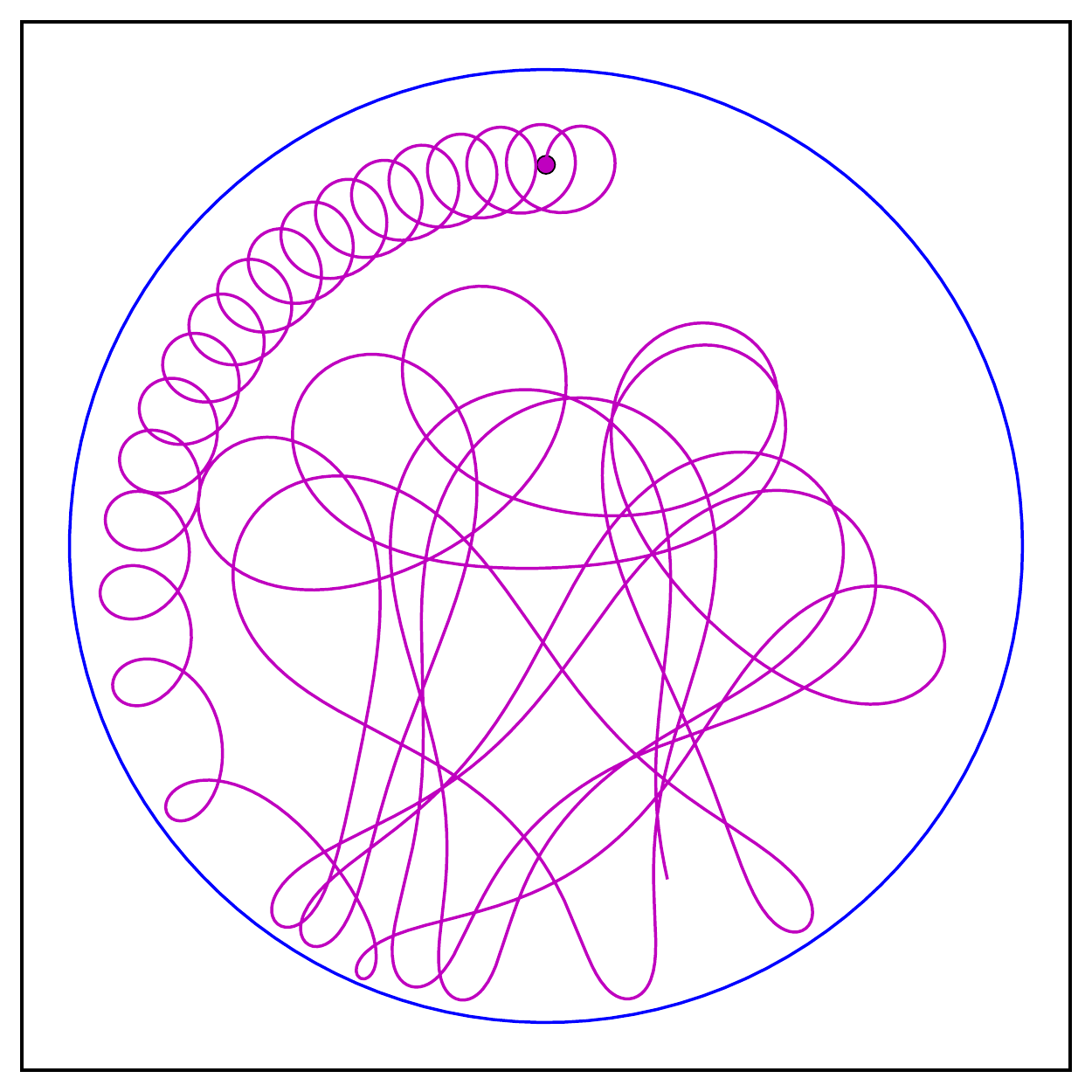}
		\caption{$B=\frac{1}{1-r}+7y+5x^2$}\label{fig:2a}		
	\end{subfigure}
	\quad
	\begin{subfigure}[t]{0.45\textwidth}
		\centering
		\includegraphics[width=\textwidth]{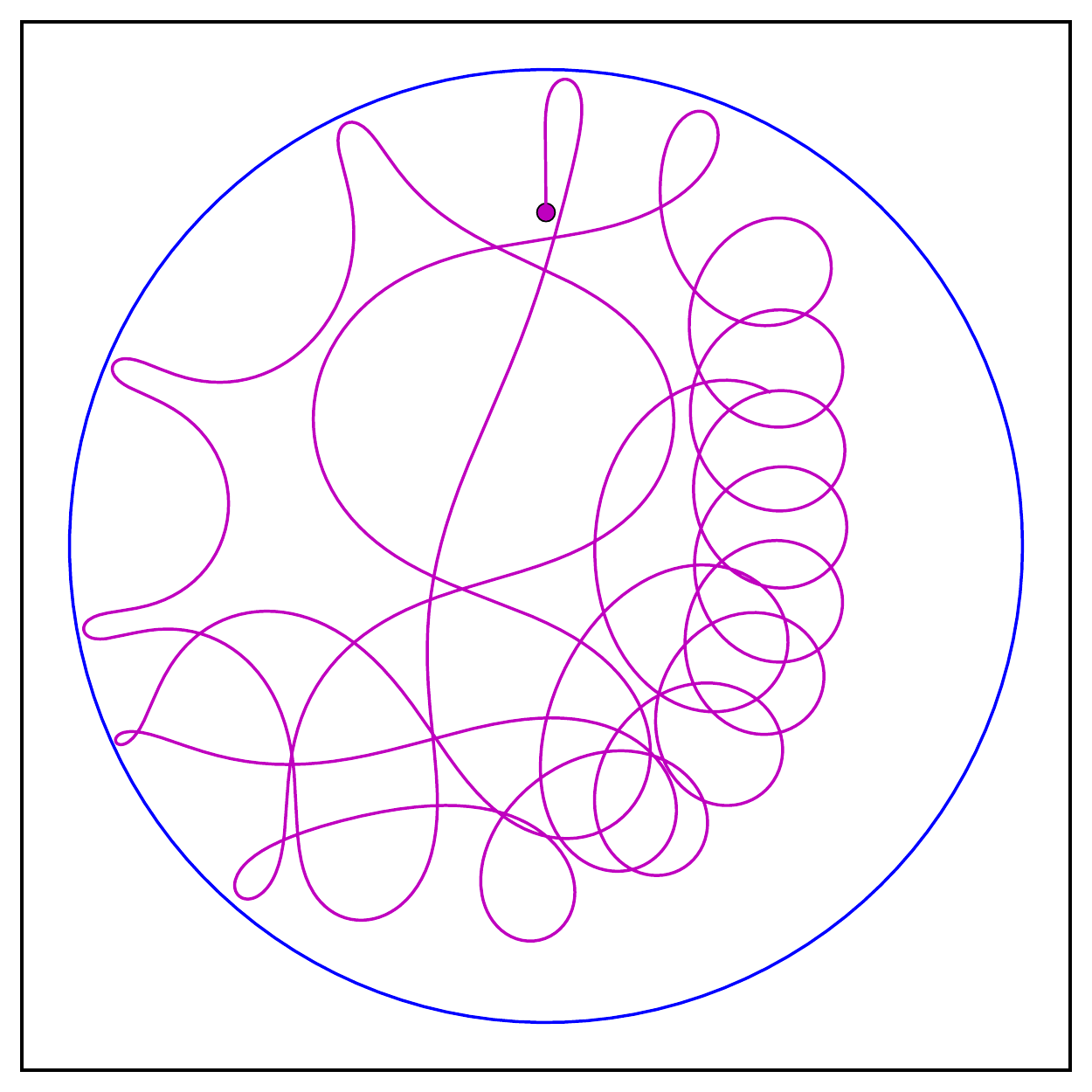}
		\caption{$B=\frac{1}{1-r}+10x-2x^2-10y^3$}\label{fig:2b}
	\end{subfigure}
	\caption{More examples in the unit disc}
	\label{fig:maggen}
\end{figure}

We experimented numerically with many examples in the unit disc. We can see in Figure \ref{fig:magsym} a couple of trajectories of charged particles in a magnetic field that depends only on $r=|q|$. The dynamics in this case is totally integrable. In figure \ref{fig:maggen} we can see a couple more examples of magnetic fields satisfying our hypothesis, notice how the particles always bounce back from the boundary.

\begin{figure}
  \centering
  \resizebox{0.75\textwidth}{!}{
  \includegraphics{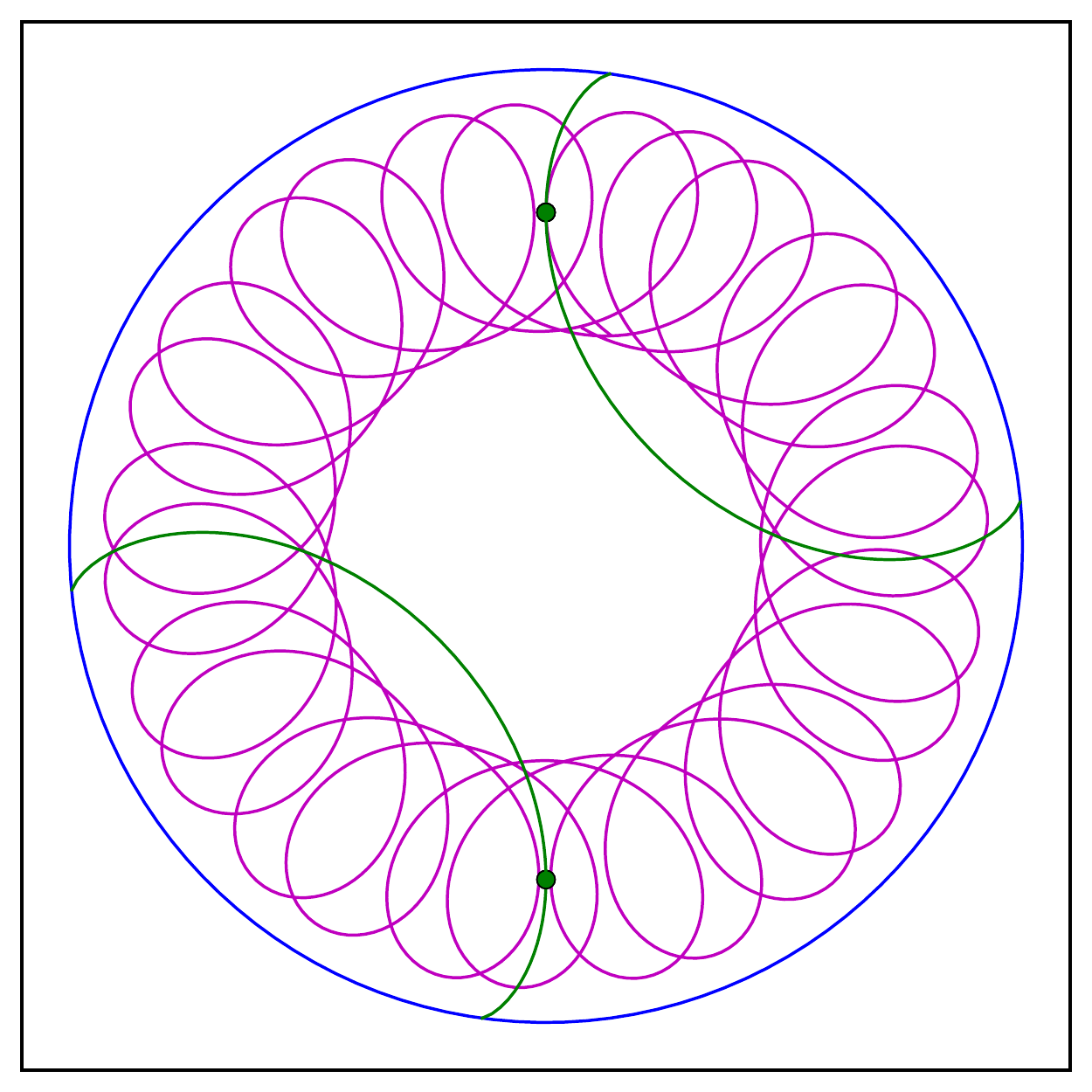}
  }
  \caption{$B(x,y)=\frac{1}{\sqrt{1-r}}$}
  \label{fig:magnogood}
\end{figure}

Finally figure \ref{fig:magnogood} shows three different trajectories of particles under a magnetic field which is integrable as $|q|\to1$, this is an example that doesn't fit our hypothesis. One of the trajectories exhibits a similar behavior to the trajectories in figure \ref{fig:magsym} and it's confined inside the disc, while the two other trajectories look like arcs with endpoints at the boundary, it seems like the magnetic field is not strong enough to push these two trajectories away from the boundary curve.

When experimenting numerically with non integrable magnetic fields, we could not find trajectories that approach the boundary as the ``arc-like'' trajectories in figure \ref{fig:magnogood}, which is a good confirmation of our result, but a relevant point is that our experiments suggest these ``arc-like'' orbits might still take infinite time to reach the boundary.

The equations were solved numerically by using the odeint integrator in the SciPy library for Python, which itself is an implementation of the LSODA integrator from the FORTRAN library ODEPACK. The figures were created with the matplotlib library for Python.

\section{Proof of our results}
\label{sec:proof}
We now use normal coordinates to find a convenient potential for the magnetic field $\mathbf{B}=B(n,s)dn\wedge ds$, let's define $\mathbf{A}=A(n,s)ds$ where
  \begin{equation}
  A(n,s) = -\int_n^N B(m,s)dm
  \end{equation}
since $\partial A/\partial n = B$ we have the desired identity $d\mathbf{A}=\mathbf{B}$. Notice that the Euclidian metric can be written in normal coordinates as $dn^2 + (1-\kappa(s)n)^2ds^2$, where $\kappa(s)$ denotes the curvature of $\mathcal{C}$, so that the Hamiltonian in these coordinates looks like:
  \begin{equation} \label{eq:FHam}
  H(q,p) = \frac{p_n^2}{2m} + \frac{(p_s-eA(n,s))^2}{2m(1-\kappa(s)n)^2}
  \end{equation}

The strategy of our proof is to show that the potential $A(n,s)$ must remain bounded along a magnetic trajectory that stays close to the boundary. This in turn implies that the particle cannot approach the boundary since the non integrability condition (\ref{eq:nintB}) in the theorem can be recast in terms of the potential as:
  \begin{equation}
  \lim_{n\to0}|A(n,s)| = \infty
  \end{equation} 

The foundations of our results lie in the following:


\begin{proposition}\label{potest}
Assume the magnetic field $\mathbf{B}$ satisfies the hypothesis of theorem \ref{main}. Assume further that $N$ is small enough so that we may define normal coordinates on $\Omega_\mathcal{C}(N)$ and such that $N<\varepsilon/K$, for some $0<\varepsilon<1$. Let $K=\sup_s|\kappa(s)|$ and $K'=\sup_s|\kappa'(s)|$. If $q(t)\subset\Omega_\mathcal{C}(N)$ is a trajectory with energy $H_0$, then:
  \begin{equation}
  |A(t)| \leq C_0 + C_1|t|
  \end{equation}
where $C_0,C_1>0$ are the following explicit constants:
  \begin{equation}
  C_0 = \frac{|p_s(0)|+\sqrt{2mH_0}(1+\varepsilon)}{e},\quad
  C_1 = \sqrt{\frac{2H_0}{m}}\cdot\frac{D_\mathcal{C}}{1-\varepsilon}+\frac{2H_0K'N}{e(1-\varepsilon)}
  \end{equation}
\end{proposition}


\textbf{Proof:} Notice that by equation (\ref{eq:FHam}), a particle with energy $H_0$ must satisfy:
  \begin{equation}
  \frac{(p_s-eA(n,s))^2}{2m(1-\kappa(s)n)^2} \leq H_0
  \end{equation}

 By our assumptions we have the following bound:
  \begin{equation}
  |p_s-eA(n,s)| \leq \sqrt{2mH_0}(1-\kappa(s)n) \leq \sqrt{2mH_0}(1+\varepsilon)
  \end{equation}
which shows that we may control $A(n,s)$ by controlling $p_s$. In order to do this, we use Hamilton's equation to see how the quantity $p_s$ evolves:  
  \begin{equation}
  \dot{p}_s = -\frac{\partial H}{\partial s} = \frac{(p_s-eA(n,s))}{m(1-\kappa(s)n)^2}\cdot e\frac{\partial A}{\partial s}
                    - \frac{(p_s-eA(n,s))^2}{m(1-\kappa(s)n)^3}\cdot\kappa'(s)n
  \end{equation}
  
Now, notice that $|\partial A/\partial s| < D_\mathcal{C}$ since:
  \begin{equation*}
  \begin{aligned}
  \left|\frac{\partial A}{\partial s}(n,s)\right| &= \left|\int_n^N \frac{\partial B}{\partial s}(m,s) dm\right| \\
                                                                 &\leq \int_n^N \left|\frac{\partial B}{\partial s}(m,s)\right| dm \\
                                                                 &\leq \int_0^N \left|\frac{\partial B}{\partial s}(m,s)\right| dm \\
                                                                 & \leq \sup_{s} \int_0^N \left|\frac{\partial B}{\partial s}(m,s)\right| dm \\
                                                                 & = D_\mathcal{C}
  \end{aligned}
  \end{equation*}
so we obtain the following bound along $q(t)$:
  \begin{equation} \label{eq:pdot}
  |\dot{p}_s| \leq \sqrt{\frac{2H_0}{m}}\cdot\frac{eD_\mathcal{C}}{1-\varepsilon} + \frac{2H_0K'N}{1-\varepsilon} = eC_1
  \end{equation}

Integrating inequality (\ref{eq:pdot}) from time $0$ to time $t$ we obtain the following:
  \begin{equation}
  |p_s(t)| \leq |p_s(0)| + eC_1|t|
  \end{equation}
  
We may then deduce:
  \begin{equation}
  \begin{aligned}
  |A(q(t))| & \leq \frac{|p_s(t) - eA(q(t)|}{e} + \frac{|p_s(t)|}{e} \\
               & \leq \frac{\sqrt{2mH_0}(1+\varepsilon)}{e} + \frac{|p_s(0)|}{e} + C_1|t|\\
               & = C_0 + C_1|t|
  \end{aligned}
  \end{equation}
which is the estimate we wanted. $\square$\\

Theorem \ref{quant} is now a simple corollary of the above proposition.

We now move on to the proof of the main theorem:\\


\textbf{Proof of Theorem \ref{main}:} Suppose we have a trajectory $q:[0,T]\to\Omega$ that reaches the boundary at time $T$ at a connected component $\mathcal{C}$. By the first condition in our theorem, i.e. equation (\ref{eq:nintB}), we must have:
  \begin{equation} \label{eq:Ablows}
  \lim_{t\to T}A(q(t)) = \infty
  \end{equation}
we may suppose further,without loss of generality, that $q(t)$ lies in $\Omega_\mathcal{C}(N)$ for all $t\in[0,T]$ for some $N$ that satisfies the hypothesis of propostion \ref{potest}. We have then:
  \begin{equation}
  |A(q(t))| \leq C_0+C_1T
  \end{equation}

Hence we conclude that $A$ is bounded along the trajectory, which contradicts (\ref{eq:Ablows}) so that no such trajectory may exist. $\square$


\section{Relations to Quantum Systems}\label{sec:quantum}

We discuss in this section an example that illustrates an interesting dichotomy between the behavior of a classical system and its quantization. To be more precise, we will describe a family of magnetic fields $\mathbf{B}_\alpha$ defined over the unit disc for which the Hamiltonian dynamics of a classical particle under the influence of $\mathbf{B}_\alpha$ is complete, but the dynamics of a quantum particle under the influence of the same field is not.

From now on denote by $\Omega$ the unit disc in $\mathbb{R}^2$
  \begin{equation}
  \Omega = \{(x,y) \in \mathbb{R}^2 \mid x^2+y^2 = r^2 < 1\}
  \end{equation}

Given a constant $\alpha>0$ consider the magnetic field:
  \begin{equation}
  \mathbf{B}_{\alpha}(x,y) = \alpha\frac{r-2}{(r-1)^2}dx\wedge dy
  \end{equation}

If we choose a potential 1-form $\mathbf{A}=A_xdx+A_ydy$ for $\mathbf{B}_\alpha$, we may define the magnetic Schr\"{o}dinger operator $\hat{H}_{\mathbf{A}}$ as the quantization of the classical Hamiltonian defined in equation (\ref{eq:Ham}), i.e.:
  \begin{equation}
  \hat{H}_{\mathbf{A}} = - (\partial_x - ieA_x)^2 - (\partial_y - ieA_y)^2
  \end{equation}
Where we set $\hbar=1$ for convenience, which can be accomplished by choosing appropriate units.

A quantum particle inside $\Omega$ is modeled by a function $\psi_0(x)$ in $L^2(\Omega)$, with $|\psi_0|_{L^2(\Omega)}=1$ and the corresponding system is governed by Schr\"{o}dinger's equation:
  \begin{equation}\label{eq:schrod}
  i\frac{\partial\psi}{\partial t}=\frac{1}{2m}\hat{H}_{\mathbf{A}}\psi
  \end{equation}
We say that the quantum system is complete or that the quantum particle is confined by the magnetic field $\mathbf{B}$ if the operator $\hat{H}_{\mathbf{A}}$ with domain $C^\infty_0(\Omega)$ of smooth functions compactly supported in $\Omega$ is essentially self-adjoint (this property does not depend on the choice of potential $\mathbf{A}$, see \cite{deVerdiere}). In this case there is a unique self-adjoint extension of $\hat{H}_{\mathbf{A}}$ and by Stone's theorem a unique strongly continuous unitary one-parameter subgroup $(U_t)_{t\in\mathbb{R}}$ so that $\psi(t,x)=U_t\psi_0(x)$ is a solution to equation (\ref{eq:schrod}).

Notice that this is a good definition of confinement. Since the operators $U_t$ are unitary we have
  \begin{equation}
  |\psi_t|_{L^2(\Omega)}=|U_t\psi_0|_{L^2(\Omega)}=|\psi_0|_{L^2(\Omega)}=1
  \end{equation}
which means that for any time $t$, the particle is observed inside of $\Omega$ with probability $1$, since $|\psi_t|^2$ is the probability density for the position operators.

In \cite{deVerdiere} Colin de Verdi\`{e}re and Truc prove that the operator $\hat{H}_\mathbf{A}$ associated to $\mathbf{B}_\alpha$ is not essentially self-adjoint in $C^\infty_0(\Omega)$ if the constant $\alpha$ satisfies $0<\alpha<\sqrt{3}/2$. So that in order to solve Schr\"{o}dinger's equation we must impose boundary conditions on $\psi$. This means that as time evolves certain wave functions will interact with the boundary, so not all quantum particles are confined to $\Omega$. In this sense we'll say the quantum dynamics is not complete

On the other hand by parametrizing the boundary by arclength through $\gamma(s)=(\cos(s),\sin(s))$, using normal coordinates $(n,s)$ we can express the magnetic field by (notice $n=1-r$):
  \begin{equation}
  \mathbf{B}_\alpha = \alpha\left(\frac{1}{n^2} - 1\right)dn\wedge ds
  \end{equation}
Which fits the hypothesis of our theorem and hence for a classical particle, the dynamics is complete.

One interesting remaining question is to prove that for magnetic fields $\mathbf{B}_\alpha$ for which the corresponding Hamiltonian operator is not essentially self-adjoint in $C^\infty_0(\Omega)$, there are solutions $\psi_t$ in $L^2(\mathbb{R}^2)$ for which the norm restricted to $\Omega$ decreases with time, so that some particles would in fact ``leak'' through the boundary walls, which would confirm a quantum tunneling property for these systems.\\


\textbf{Acknowledgements}\\

I would like to thank Yves Colin de Verdi\`{e}re and Fran\c{c}oise Truc for having the generosity of sharing this interesting problem in their work. I'd also like to thank Richard Montgomery for his invaluable guidance and encouragement during this project.


\end{document}